\documentclass[12pt]{article}

\usepackage{amsmath}
\usepackage{calc}

\makeatletter
\newlength{\raise@amount@commafrac}
\newcommand{\textraised}[1]{%
  \setlength{\raise@amount@commafrac}%
      {\heightof{$\vcenter{\hbox{$.$}}$}-\heightof{#1}}%
  \raisebox{\raise@amount@commafrac}{#1}%
}
\makeatother

\date{}

\begin{document}

\centerline{}

\centerline{\Large{\bf On the transcendence of some operations}}
\centerline{\Large{\bf of infinite series}}

\centerline{}


\newcommand{\mvec}[1]{\mbox{\bfseries\itshape #1}}

\centerline{\bf {Sarra Ahallal, Fedoua Sghiouer, Ali Kacha,}}

\centerline{}

\centerline{Department of Mathematics, Faculty of Science,}

\centerline{Ibn Tofail University, 14 000 Kenitra, Morocco  }

\centerline{e-mail: sarra.ahallal@uit.ac.ma}

\centerline{e-mail: fedoua.sghiouer@uit.ac.ma} 

\centerline{e-mail: ali.kacha@uit.ac.ma}

\newtheorem{Theorem}{\quad Theorem}[section]

\newtheorem{Definition}[Theorem]{\quad Definition}

\newtheorem{Corollary}[Theorem]{\quad Corollary}

\newtheorem{Lemma}[Theorem]{\quad Lemma}

\newtheorem{Example}[Theorem]{\quad Example}

 \begin{center}
     \textbf{\large{Abstract}}
 \end{center}

In the present paper and as an application of Roth's theorem concerning the rational approximation of algebraic numbers, we give a sufficient condition that will assure us that a sum, product and quotient of some series of positive rational terms are transcendental numbers.
We recall that all the infinite series that we are going to treat are Liouville numbers. At the end this article, we establish an approximation measure of these numbers.
 \\

\hspace*{-0.6cm}\textbf{Keywords:} Infinite series, transcendental number, approximation measure.

\section{Introduction and preliminaries }

The theory of transcendental numbers has a long history. We know since J. Liouville in 1844 that the very rapidly converging sequences of rational numbers provide examples of transcendental numbers. So, in his famous paper [8], Liouville showed that a real number admitting very good rational approximation can not be algebraic, then he explicitly constructed the first examples of transcendental numbers. In [1], the authors give a family of real numbers that are transcendental.\\

There are a number of sufficient conditions known within the literature for an infinite series, $ \sum_{n = 1}^{\infty} 1/a_n$, of positive rational numbers to converge to an irrational number, see [3, 10].These conditions, which are quite varied in form, share one common feature, namely, they all require rapid growth of the sequence ${ a_n }$ to deduce irrationality of the series. As an illustration consider the following results of J. Sandor which have been taken from [12] and [13].\\

From this direction, the transcendence of some infinite series has been studied by several authors such as M.A. Nyblom [9], J. Hančl and J. Štěpnička [5]. we also note that the transcendence of some power series with rational coefficients is given by some authors such as J. P. Allouche [3] and G. K. Gözer [4].
\\
We recall that in $[5],$ we have proved the transcendence of some series of positive rational terms. \\
\\
\noindent In the present paper, the first aim is to give a sufficient condition that will assure us that a sum, product and quotient of some series of positive rational terms are transcendental numbers.
\\
 \noindent Let $g_1, g_2$ be two distinct integers $\geq 2,$ such that $g_1>g_2$ and $\theta_1,$ $ \theta_2$ two infinite series
 which are defined by
$$
\theta_1=\sum_{n=1}^{+\infty}g_1^{-a_n}, \ \theta_2=\sum_{n=1}^{+\infty}g_2^{-a_n}
$$
\noindent where $a_{n} \geq 1$, $b_{n} \geq 1 $ are integers for all $n\geq 1.$
\\
\\
\noindent The second main result of this article is to establish an
approximation measure of a real number $\theta_1+\theta_2.$
\\
\noindent In order to prove the transcendence of the infinite series, we will use the Roth's Theorem.

{\bigskip}

\noindent {\bf  Roth's Theorem. \ }
 {\it Let $ \alpha $ be a real number, $ \delta $ a real number $> 2,$ if there exists an infinity rational numbers $ \frac{p}{q} $ with $ gcd(p, q) = 1 $ such that
 $$ \left| \alpha - \frac{p}{q} \right| < \frac{1}{q^{ \delta} },
  $$

\noindent then $ \alpha $ is a transcendental number.
}

{\bigskip}

\section{Main results}

{\bigskip}

\noindent{\bf 2.1 Transcendence of some series }

{\bigskip}

\noindent{\bf Theorem 1. \ } {\it Let $g_1,g_2$ be two distinct integers $\geq 2,$ such that $g_1 > g_2,$ $\beta$ a real number $>0.$ With the same notations as above,
let
$$\theta_1=\sum_{n=1}^{+\infty}\frac{1}{g_1^{a_n}}, \text{ } \theta_2=\sum_{n=1}^{+\infty}\frac{1}{g_2^{a_n}},
$$
\noindent where $a_{n+1}=a_n^{1+\beta}$ for all $n \geq 1 $  and $a_1 \geq 2.$
\\
\noindent Then, the real numbers $\theta_1+\theta_2,$ $\theta_1-\theta_2,$ $\theta_1 \times \theta_2 $ and $ \theta_1/\theta_2$ are transcendental numbers.}

{\bigskip}

\noindent{\bf Remark $1.$ }  From the definition of $\theta_j,$ it is clear that these series are convergent.
\\
\\
\noindent{\bf Proof of Theorem $1.$ \ } Let
 $$
 \theta_{n,1}=\sum_{k=1}^{n} \frac{1}{g_1^{a_k}} =\frac{p_{n,1}}{q_{n,1}}, \\
 \\
  \theta_{n,2}=\sum_{k=1}^{n} \frac{1}{g_2^{a_k}} =\frac{p_{n,2}}{q_{n,2}},
 $$

 \noindent where $\gcd(p_{n,1}, q_{n,1})= \gcd(p_{n,2}, q_{n,2})=1.$  So, we can see that $ q_{n,1}=g_1^{a_n}$ and  $ q_{n,2}=g_2^{a_n}.$

 {\bigskip}

 \noindent (i) Transcendence of  $(\theta_1+\theta_2).$  \\
 \\
 We have
 $$
 |\theta_1+\theta_2-\theta_{n,1}-\theta_{n,2}| \leq |\theta_1-\theta_{n,1}| +|\theta_2-\theta_{n,2}|.
 $$

 \noindent Therefore, we get
 $$
 |\theta_1-\theta_{n,1}|=\frac{1}{g_1^{a_{n+1}}} \left( 1+\sum_{k=n+1}^{+\infty} \frac{1}{g_1^{a_{k+1}-a_{n+1}}}\right).
 $$
 \noindent which gives
 $$
  \frac{1}{g_1^{a_{n+1}}}  < |\theta_1-\theta_{n,1}|<\frac{2}{g_1^{a_{n+1}}},
  $$
 \noindent and
 $$
 \frac{1}{g_2^{a_{n+1}}}  < |\theta_2-\theta_{n,2}|<\frac{2}{g_2^{a_{n+1}}}.
 $$
 \noindent From the previous inequalities and the hypothesis $g_1 > g_2,$ we deduce

 \begin{equation}
 |\theta_1+\theta_2-\theta_{n,1}-\theta_{n,2}|  < \frac{2}{g_1^{a_{n+1}}}+ \frac{2}{g_2^{a_{n+1}}}
 \\
 \\
                                                 <   \frac{4}{g_2^{a_{n+1}}}.
 \end{equation}

\noindent We can see that
$$
\lim_{n \rightarrow +\infty} \frac{g_2^{a_{n+1}}}{(g_1g_2)^{a_{n}}}=+\infty.
$$
\noindent The above inequality is true since one has
$$
\frac{a_{n+1} \ln g_2}{a_n \ln (g_1g_2)}=a_n^{\beta}\frac{\ln g_2}{\ln (g_1g_2)} \geq 2^{\beta (1+\delta)^{n-1}} \frac{\ln g_2}{\ln (g_1g_2)}.
$$

\noindent and $ \lim_{n \rightarrow +\infty} 2^{\beta (1+\delta)^{n-1}}\frac{\ln g_2}{\ln (g_1g_2)}= +\infty.
$  \\
\\
\noindent Hence, for all positive integer $d >2,$ $\exists n_0=n_0(d)$ such that for all $ n \geq n_0,$ we have
$$
\ln (g_2^{a_{n+1}}) > d \ln (g_1g_2)^{a_n}= \ln {((g_1g_2)^{a_n})}^d.
$$
\noindent Which yields,
\begin{equation}
g_2^{a_{n+1}} > (g_1^{a_n}g_2^{a_n})^d
\end{equation}

\noindent  and  the inequality $(1)$ becomes
$$
|\theta_1+\theta_2-(\frac{p_{n,1}}{q_{n,1}}+\frac{p_{n,2}}{q_{n,2}})|< \frac{4}{(q_{n,1}q_{n,2})^d}.
$$

\noindent Since, $ d >2,$ by Roth's Theorem, we deduce that $\theta_1+\theta_2$ is a transcendental number.
\\
\noindent Similarely, it can easily be proven that $\theta_1-\theta_2$ is also a transcendental number.

{\bigskip}

\noindent (ii) {  Transcendence of  $\theta_1\times \theta_2.$
\\
\noindent One has
\begin{align*}
|\theta_1 \times \theta_2-\theta_{n,1} \times \theta_{n,2}| &= |(\theta_1-\theta_{n,1}) \theta_{n,2}+ (\theta_{n,2}- \theta_2)\theta_{n,1}|
                                                            \\
                                                            &\leq  |\theta_1-\theta_{n,1}| \theta_2+(1+\theta_1) |\theta_2-\theta_{n,2}|
                                                            \\
                                                            &<\frac{2 \theta_2}{g_1^{a_{n+1}}}+\frac{2(1+\theta_1)}{g_2^{a_{n+1}}}
\end{align*}

\noindent for all $n$ sufficiently large. Since $g_1 > g_2$ and using the relationship (2), we obtain
$$
|\theta_1 \times \theta_2-\theta_{n,1} \times \theta_{n,2}| < \frac{2(1+\theta_1+\theta_2)}{g_2^{a_{n+1}}} <  \frac{1}{(q_{n,1}q_{n,2})^d},
$$
\noindent for all $n$ sufficiently large. Therefore, $\theta_1 \times \theta_2$ is transcendental by Roth's Theorem.

{\bigskip}

\noindent (iii) { Transcendence of  $ (\theta_1/ \theta_2).$

\begin{align*}
\left|\frac{\theta_1}{ \theta_2}-\frac{\theta_{n,1}}{\theta_{n,2}}\right| &= |(\theta_1-\theta_{n,1}) \theta_2+ (\theta_{n,2}- \theta_2)\theta_{n,1}|
                                                            \\
                                                            &=  \frac{|\theta_{n,2}(\theta_1-\theta_{n,1})+\theta_{n,1}(\theta_{n,2}-\theta_{2})|}{\theta_2\theta_{n,2}}
 \end{align*}
 \noindent As $\theta_{n,1}={p_{n,1}}/{q_{n,1}} \leq 1+\theta_1< 2 $ and $ \theta_{n,2}= {p_{n,2}}/{q_{n,2}} > {1}/{g_2^{a_1}},$ for
 $n \geq 2,$ we obtain

\begin{align*}
\left|\frac{\theta_1}{ \theta_2}-\frac{\theta_{n,1}}{\theta_{n,2}}\right| & \leq \frac{|\theta_1-\theta_{n,1}|}{\theta_2}+\frac{2g_2^{a_1}}{\theta_2} |\theta_{n,2}- \theta_2|
                                                                           \\
                                                                           &<\frac{2 }{\theta_2 g_1^{a_{n+1}}}+\frac{4g_2^{a_1}}{\theta_2 g_2^{a_{n+1}}}
                                                                           \\
                                                                           & < \frac{2+4g_2^{a_1}}{\theta_2}\frac{1}{g_2^{a_{n+1}}}.
\end{align*}

\noindent Since $p_{n,2} < (1+\theta_2) q_{n,2} < 2 q_{n,2}$ and using $(2),$ for all $n$ sufficiently large,
\noindent we find
\begin{equation*}
\left|\frac{\theta_1}{ \theta_2}-\frac{\theta_{n,1}}{\theta_{n,2}}\right|  \leq \frac{4}{(q_{n,1}q_{n,2})^d} < \frac{4(1+\theta_2)}{(q_{n,1}p_{n,2})^d}.
\end{equation*}

\noindent Therefore ${\theta_1}/{\theta_2}$ is a transcendental number.
{\bigskip}

\noindent{\bf  2.2. Approximation measure of a power series $(\theta_1+\theta_2 )$}

{\bigskip}

\noindent In this subsection, we give the second main result of
this article. \noindent Throughout this section, we adopt the following
notation
$\theta=\theta_1 + \theta_2 $ and $\theta_n = \theta_{1,n} + \theta_{2,n}= p_n/q_n$ where
$\gcd (p_n,q_n)=1.$

{\bigskip}

\noindent {\it {\bf Theorem 2.\ }
\noindent  Let $\xi$ be an algebraic number of degree $d\geq 2$ and height $%
H.$  Let $\alpha > H$ and $k >1 $ be two real numbers such that for all $n\geq 1, $ we have
\begin{equation*}
a_{n}^{\alpha } \leq a_{n+1}<a_{n}^{k\alpha },
\text{ \ for all \ } n\geq 1.%
\end{equation*}
\noindent Then we get,
\begin{equation*}
\left\vert \theta-\xi \right\vert > \frac{1}{(2Hd^2)^{1+4d }}.
\end{equation*}
}
 {\bigskip}

 \noindent {\bf Remark $2.$ } This approximation measure obtained of $\theta_1+\theta_2$ is the same of $\theta_1-\theta_2,$ $\theta_1\times \theta_2$ and $ \theta_1/\theta_2.$

 {\bigskip}

\noindent {\bf Proof of Theorem 2. \ } Let $\xi $ be an algebraic
number of degree $d\geq 2$ and height $H >\alpha.$
From Theorem $1,$ part $i)$ we deduce that $\theta$ is transcendental.
\\
\noindent To obtain an approximation measure of $\theta_1+\theta_2$, it is
sufficient to minimize $ \left\vert \theta-\xi \right\vert $ by
function of $H$ and $d.$

{\bigskip}

\noindent We have
$$
\vert \theta-\xi|=  \vert \theta -\theta_n + \theta_n-\xi \vert \leq
\vert \theta -\theta_n \vert +\vert \xi-\theta_n \vert.
$$
Which becomes,
$$
\vert \theta-\xi \vert  \geq \vert \theta - \xi \vert -\vert \theta -\theta_n \vert.
$$
\noindent  According to Liouville's Theorem, it is not hard to see
that
\begin{equation}
\left\vert \xi -\theta_n%
\right\vert \geq \frac{1}{Hd^{2}q_{n}^{d }}.
\end{equation}

\noindent It follows from section $1$ that
\begin{equation*}
\left\vert \theta -\theta_{n}\right\vert
<\frac{1}{g_2^{a_{n+1}}}.
\end{equation*}

\noindent so, we obtain
\begin{equation*}
\vert \xi -\theta \vert >  \frac{1}{Hd^{2}q_{n}^{d }}-\frac{1}{g_2^{a_{n+1}}}.
\end{equation*}

\noindent In order to have
\begin{equation*}
\vert \xi -\theta \vert > \frac{1}{2Hd^{2}q_{n}^{d }},
\end{equation*}
it suffices to have
\begin{equation*}
 \frac{1}{2Hd^{2}q_{n}^{d }}> \frac{1}{g_2{^a_{n+1}}}.
\end{equation*}

\noindent Since $q_n=(g_1g_2)^{a_n},$ the relationship $(3)$ gives
\begin{equation}
 \frac{1}{2Hd^{2}q_{n}^{d }}> \frac{1}{(g_1g_2)^{a_{n+1}}}.
\end{equation}

\noindent Therefore $(4)$ is equivalent to

\begin{equation}
 (g_1g_2)^{a_{n+1}} > 2Hd^2(g_1g_2)^{da_n }.
\end{equation}

\noindent To realize the inequality $(5),$ it suffices that the
integer $n$ satisfies

\begin{equation}
 \left\{
    \begin{array}{ll}
     (g_1g_2)^{a_{n+1}/2} > (g_1g_2)^{da_n }\\
     \\
     (g_1g_2)^{a_{n+1}/2} > 2Hd^2.
    \end{array}
\right.
\end{equation}

\noindent The first inequality of $(6)$ is easily obtained because
$a_{n+1} > a_n^{\alpha} > 2da_n.$
\\
\noindent For the second one, let $n_1$ be the smallest integer $n$ such that
\begin{equation}
 (g_1g_2)^{a_{n}/2} < 2Hd^2 < (g_1g_2)^{a_{n+1}/2 }.
 \end{equation}

\noindent This yields that
\begin{equation*}
\vert A-\xi \vert > \frac{1}{\left(
2Hd^{2}\right) ^{1+4d}}.
\end{equation*}

\noindent Which proves the result of Theorem $2.$

{\bigskip}

\noindent {\bf Example. \ } Let $$ \left\{
    \begin{array}{ll}
     g_1 = 2, g_2 = 3, a_1=2
     \\
        a_{n+1} = a_n^2, n \geq 1,
         \\
        \alpha = 4, k = 2.
    \end{array}
\right. $$
Let $\xi$ be an algebraic number of degree $d=3$ and height $H.$ 
 By applying Theorem $2,$ an approximation measure of $ (\theta_1+ \theta_2 )$ is given by

    $$
    \left|\theta_1+\theta_2-\xi \right| > \dfrac{1}{(18H)^{13}}.
    $$

{\bigskip}

\end{document}